\documentclass[graybox]{svmult}
\pdfoutput=1
\usepackage{type1cm}   
\usepackage{makeidx} 
\usepackage{graphicx} 
\usepackage{multicol}  
\usepackage[bottom]{footmisc}
\usepackage{newtxtext}          
\usepackage{newtxmath}
\usepackage{subcaption}
\usepackage{caption}
\usepackage{nicefrac}
\usepackage{hyperref}
\usepackage[numbers,sort&compress]{natbib}

\allowdisplaybreaks[4]

\begin{document}

\title*{Simultaneous approximation terms for elastic wave equations on nonuniform grids}
\author{Longfei Gao and David Keyes}
\institute{Longfei Gao \& David Keyes \at Division of Computer, Electrical and Mathematical Sciences and Engineering, 
       \ \\ King Abdullah University of Science and Technology, Thuwal 23955-6900, Saudi Arabia
       \ \\ e-mail: longfei.gao@kaust.edu.sa \& david.keyes@kaust.edu.sa }

\maketitle

\vspace{-7.em}
\abstract{We consider the finite difference discretization of isotropic elastic wave equations on nonuniform grids. The intended applications are seismic studies, where heterogeneity of the earth media can lead to severe oversampling for simulations on uniform grids. To address this issue, we demonstrate how to properly couple two non-overlapping neighboring subdomains that are discretized uniformly, but with different grid spacings. Specifically, a numerical procedure is presented to impose the interface conditions weakly through carefully designed penalty terms, such that the overall semi-discretization conserves a discrete energy resembling the continuous energy possessed by the elastic wave system.}


\section{Introduction}
Numerical simulation of wave phenomena is routinely used in seismic studies, where simulated wave signals are compared against experimental ones to infer subterranean information. 
Various systems can be used to model wave propagation in earth media.
Here, we consider the system of isotropic elastic wave equations described in Section \ref{Section_problem_description}.
Various numerical methods can be applied to discretize such system, among which the finite difference methods (FDMs) are still very popular, particularly for seismic exploration applications, due to their simplicity and efficiency.

However, when discretized on uniform grids, heterogeneity of the earth media will lead to oversampling in both space and time, undermining the efficiency of FDMs.
Specifically, since spatial grid spacing is usually decided on a point-per-wavelength basis for wave simulation, uniform grid discretization will lead to oversampling in space for regions with higher wave-speeds.
On the other hand, temporal step length is usually restricted by the Courant-Friedrichs-Lewy (CFL) stability condition for wave simulation, which will lead to oversampling in time for regions with lower wave-speeds. 

For earth media, the contrast between smallest and largest wave-speeds can be as high as fifty (see \cite[~p.240]{thierry1987acoustics}), which entails significant oversampling for discretizations on uniform grids.
Moreover, wave-speeds in earth media tend to increase with depth due to sedimentation and consolidation.
These observations motivate us to consider the grid configuration illustrated in Figure \ref{Grid_Configuration_Interface}, where two uniform grid regions are separated by a horizontal interface. 
The staggered grids approach (see, e.g., \cite{yee1966numerical}) is considered here, which uses different subgrids to discretize different solution variables. 
In Figure \ref{Grid_Configuration_Interface}, the ratio of the grid spacings in the two regions is two. However, other ratios, not necessarily integers, can also be addressed with the methodology presented here. 
Furthermore, multiple grid layers can be combined together in a cascading manner to account for larger wave-speed contrasts.

In this work, we recap one of the earliest motivations of domain decomposition methods in demonstrating how to combine the two regions illustrated in Figure \ref{Grid_Configuration_Interface} without numerical artifacts.
Specifically, we adopt the summation by parts (SBP) - simultaneous approximation terms (SAT) approach, which utilizes discrete energy analysis to guide the discretization.
The overall semi-discretization is shown to be discretely energy conserving, preserving the same property in the continuous elastic wave system.
The concept of SBP operators dates back to \cite{kreiss1974finite} while the technique of SAT was introduced in \cite{carpenter1994time}. 
The two recent review papers \cite{svard2014review} and \cite{fernandez2014review} provide comprehensive coverage of their developments.  

In the following, we describe the abstracted mathematical problem in Section \ref{Section_problem_description}, present the interface treatment in Section \ref{section_methodology}, provide numerical examples in Section \ref{section_examples}, and summarize in Section \ref{section_summary}.

\section{Problem Description} \label{Section_problem_description}
We consider the 2D isotropic elastic wave equations posed as the following first-order dynamical system written in terms of velocity and stress:
\begin{equation}
\label{Elastic_wave_equation_PDE_2D_implementation}
\small
\left\{
\begin{array}{rcl}
\displaystyle \frac{\partial v_x}{\partial t} &=& \displaystyle \frac{1}{\rho} \left( \frac{\partial \sigma_{xx}}{\partial x} + \frac{\partial \sigma_{xy}}{\partial y} \right); \\[2ex]
\displaystyle \frac{\partial v_y}{\partial t} &=& \displaystyle \frac{1}{\rho} \left( \frac{\partial \sigma_{xy}}{\partial x} + \frac{\partial \sigma_{yy}}{\partial y} \right); \\[2ex]
\displaystyle \frac{\partial \sigma_{xx}}{\partial t} &=& \displaystyle \left(\lambda + 2\mu\right) \frac{\partial v_x}{\partial x} + \lambda \frac{\partial v_y}{\partial y} + \mathcal{S}; \\[2ex]
\displaystyle \frac{\partial \sigma_{xy}}{\partial t} &=& \displaystyle \mu \frac{\partial v_y}{\partial x} + \mu \frac{\partial v_x}{\partial y}; \\[2ex]
\displaystyle \frac{\partial \sigma_{yy}}{\partial t} &=& \displaystyle \lambda \frac{\partial v_x}{\partial x} + \left(\lambda + 2\mu\right) \frac{\partial v_y}{\partial y} + \mathcal{S},
\end{array}
\right.
\end{equation}
where $v_x$ and $v_y$ are particle velocities; $\sigma_{xx}$, $\sigma_{xy}$ and $\sigma_{yy}$ are stress components; $\rho$, $\lambda$, $\mu$ are density, first and second Lam\'{e} parameters that characterize the media\footnotemark[1]; and $\mathcal S$ is the source term that drives the wave propagation. 
For simplicity, the source term $\mathcal S$ is omitted in the upcoming discussion.
All solution variables and their derivatives are assumed to be zero at the initial time.
We consider periodic boundary condition for left and right boundaries and free-surface boundary condition for top and bottom boundaries.
\footnotetext[1]{Lam\'{e} parameters $\lambda$ and $\mu$ are related with the compressional and shear wave-speeds $c_p$ and $c_s$ through $\lambda = \rho ( c_p^2 - 2 c_s^2 )$ and $\mu = \rho c_s^2$. For the numerical examples appearing later in Section \ref{section_examples}, $c_p$ and $c_s$ are prescribed instead of $\lambda$ and $\mu$.}

The above system is equivalent to system \eqref{Elastic_wave_equation_PDE_2D_analysis}, which is more convenient for energy analysis and formula derivation. 
In \eqref{Elastic_wave_equation_PDE_2D_analysis}, the Einstein summation convention applies to subscript indices $k$ and $l$. 
Coefficients $s_{xxkl}$, $s_{xykl}$ and $s_{yykl}$ are components of the compliance tensor, which can be expressed in terms of $\lambda$ and $\mu$. 
However, their exact expressions are not needed for the upcoming discussion. 
As explained later, system \eqref{Elastic_wave_equation_PDE_2D_implementation} is still the one used for implementation. 
\begin{equation}
\label{Elastic_wave_equation_PDE_2D_analysis}
\small
\left\{
\begin{array}{rcl}
\displaystyle \rho \frac{\partial v_x}{\partial t} &=& \displaystyle \frac{\partial \sigma_{xx}}{\partial x} + \frac{\partial \sigma_{xy}}{\partial y}; \\[2ex]
\displaystyle \rho \frac{\partial v_y}{\partial t} &=& \displaystyle \frac{\partial \sigma_{xy}}{\partial x} + \frac{\partial \sigma_{yy}}{\partial y}; \\[2ex]
\displaystyle s_{xxkl} \frac{\partial \sigma_{kl}}{\partial t} &=& \displaystyle \frac{\partial v_x}{\partial x}; \\[2ex]
\displaystyle s_{xykl} \frac{\partial \sigma_{kl}}{\partial t} &=& \displaystyle \frac{1}{2} \left( \frac{\partial v_y}{\partial x} + \frac{\partial v_x}{\partial y} \right); \\[2ex]
\displaystyle s_{yykl} \frac{\partial \sigma_{kl}}{\partial t} &=& \displaystyle \frac{\partial v_y}{\partial y}. 
\end{array}
\right.
\end{equation}

\begin{figure}[h]
\vspace{-0.75em}
\captionsetup{justification=centering}
\captionsetup{width=1\textwidth, font=small,labelfont=small}
\centering\includegraphics[scale=0.0725]{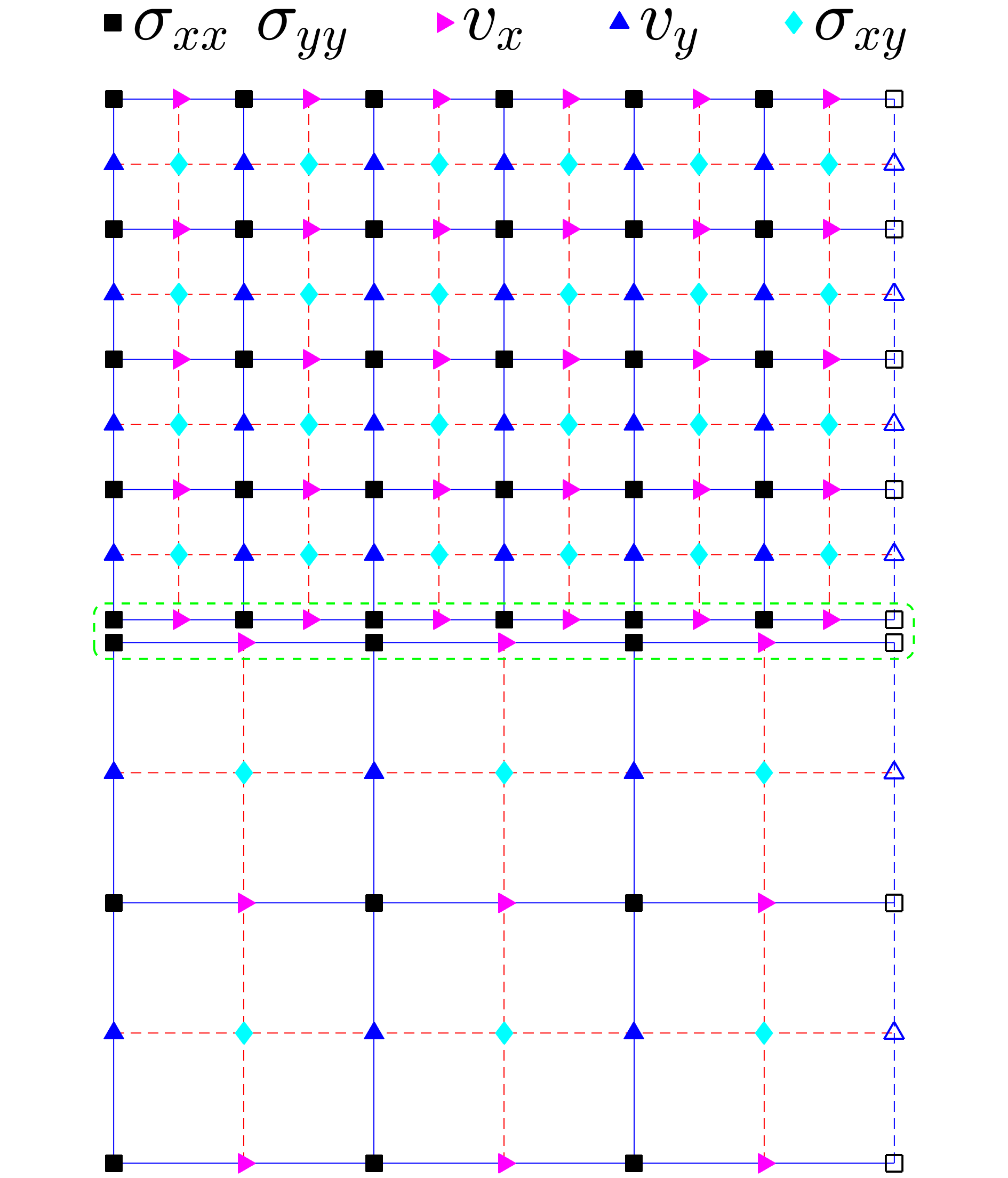}
\vspace{-0.25em}
\caption{Illustration of the grid configuration.}
\label{Grid_Configuration_Interface}
\end{figure}

The staggered grids illustrated in Figure \ref{Grid_Configuration_Interface} are used to discretize the above systems, where two uniform grid regions are separated by a horizontal interface, with a contrast ratio $1 \!\! : \!\! 2$ in grid spacing. Both regions include the interface and have grid points defined on it. 

\section{Methodology} \label{section_methodology}
In this section, we demonstrate how to couple the discretizations of system \eqref{Elastic_wave_equation_PDE_2D_implementation} on the two uniform grid regions illustrated in Figure \ref{Grid_Configuration_Interface} using the SBP-SAT approach.
A similar work has been presented in \cite{gao2019sbp} for acoustic wave equations. We will follow the methodology and terminology developed therein. 
Other related works include \cite{kozdon2016stable, wang2016high, o2017energy}.

The continuous energy associated with system \eqref{Elastic_wave_equation_PDE_2D_analysis}, and system \eqref{Elastic_wave_equation_PDE_2D_implementation} by equivalence, can be expressed as: 
\begin{equation}
\label{Continuous_energy}
\small
e = \int_\Omega \tfrac{1}{2} \rho v_i v_i d_\Omega + \int_\Omega \tfrac{1}{2} \sigma_{ij} s_{ijkl} \sigma_{kl} d_\Omega\ ,
\end{equation} 
where $\Omega$ is the simulation domain and the Einstein summation convention applies to subscript indices $i$, $j$, $k$ and $l$.
The two integrals of \eqref{Continuous_energy} correspond to the kinetic and potential parts of the continuous energy, respectively.
Differentiating $e$ with respect to time $t$ and substituting the equations in \eqref{Elastic_wave_equation_PDE_2D_analysis}, it can be shown that 
\begin{equation}
\label{Continuous_energy_conservation}
\small
\frac{d e}{d t} = \int_{\partial \Omega} v_i \sigma_{ij} n_j d_{\partial \Omega},
\end{equation} 
where $\partial \Omega$ is the boundary of the simulation domain.
For the free-surface boundary condition, i.e., $\sigma_{ij} n_j = 0$, and periodic boundary condition considered in this work, we have $\tfrac{d e}{d t} = 0$, i.e., system \eqref{Elastic_wave_equation_PDE_2D_analysis} conserves energy $e$.

Spatially discretizing \eqref{Elastic_wave_equation_PDE_2D_analysis} with finite difference methods on a uniform grid leads to the following semi-discretized system:
\begin{equation}
\label{Elastic_wave_equation_Discretization_2D_analysis}
\small
\left\{
\begin{array}{rcl}
\displaystyle \mathcal A^{V_x} \boldsymbol{\rho}^{V_x} \frac{d V_x}{d t} &=& \displaystyle \mathcal A^{V_x} \mathcal D^{\Sigma_{xx}}_x \Sigma_{xx} + \mathcal A^{V_x} \mathcal D^{\Sigma_{xy}}_y \Sigma_{xy}; \\[2ex]
\displaystyle \mathcal A^{V_y} \boldsymbol{\rho}^{V_y} \frac{d V_y}{d t} &=& \displaystyle \mathcal A^{V_y} \mathcal D^{\Sigma_{xy}}_x \Sigma_{xy} + \mathcal A^{V_y} \mathcal D^{\Sigma_{yy}}_y \Sigma_{yy}; \\[2ex]
\displaystyle \mathcal A^{\Sigma_{xx}} S^{\Sigma_{kl}}_{xxkl} \frac{d \Sigma_{kl}}{d t} &=& \displaystyle \mathcal A^{\Sigma_{xx}} \mathcal D^{V_x}_x V_x; \\[2ex]
\displaystyle \mathcal A^{\Sigma_{xy}} S^{\Sigma_{kl}}_{xykl} \frac{d \Sigma_{kl}}{d t} &=& \displaystyle \tfrac{1}{2} \mathcal A^{\Sigma_{xy}} \left( \mathcal D^{V_y}_x V_y + \mathcal D^{V_x}_y V_x \right); \\[2ex]
\displaystyle \mathcal A^{\Sigma_{yy}} S^{\Sigma_{kl}}_{yykl} \frac{d \Sigma_{kl}}{d t} &=& \displaystyle \mathcal A^{\Sigma_{yy}} \mathcal D^{V_y}_y V_y,
\end{array}
\right.
\end{equation}
where the index summation convention applies only to $k$ and $l$ in the subscripts, but not to those appearing in the superscripts. 
Superscript such as $^{V_x}$ indicates which grid the underlying quantity or operator is associated with.
In \eqref{Elastic_wave_equation_Discretization_2D_analysis}, $\mathcal D$ symbolizes a finite difference matrix, while $\mathcal A$ symbolizes a diagonal norm matrix whose diagonal component loosely representing the area that the corresponding grid point occupies. 
From the implementation perspective, the norm matrices in \eqref{Elastic_wave_equation_Discretization_2D_analysis} are redundant, but they will play an important role in deriving the proper interface treatment.
These 2D finite difference matrices and norm matrices are constructed from their 1D counterparts via tensor product. 
Specifically, 
\begin{equation}
\label{Tensor_product_norm_matrices}
\small
\begin{array}{ll}
\mathcal A^{V_x} = \mathcal A^M_x \otimes \mathcal A^N_y,\enskip & 
\mathcal A^{V_y} = \mathcal A^N_x \otimes \mathcal A^M_y, \\[0.75ex]
\mathcal A^{\Sigma_{xy}} = \mathcal A^M_x \otimes \mathcal A^M_y,\enskip &
\mathcal A^{\Sigma_{xx}} = \mathcal A^{\Sigma_{yy}} = \mathcal A^N_x \otimes \mathcal A^N_y,
\end{array}
\end{equation}
and
\begin{equation}
\label{Tensor_product_difference_matrices}
\small
\begin{array}{c}
\mathcal D_{x}^{V_x} \! = \mathcal D^M_x \otimes \mathcal I^N_y, \enskip 
\mathcal D_{x}^{V_y} \! = \mathcal D^N_x \otimes \mathcal I^M_y, \enskip 
\mathcal D_{x}^{\Sigma_{xy}} \! = \mathcal D^M_x \otimes \mathcal I^M_y, \enskip
\mathcal D_{x}^{\Sigma_{xx}} \! = \mathcal D^N_x \otimes \mathcal I^N_y,
\\[0.75ex]
\mathcal D_{y}^{V_x} \! = \mathcal I^M_x \otimes \mathcal D^N_y, \enskip 
\mathcal D_{y}^{V_y} \! = \mathcal I^N_x \otimes \mathcal D^M_y, \enskip
\mathcal D_{y}^{\Sigma_{xy}} \! = \mathcal I^M_x \otimes \mathcal D^M_y, \enskip
\mathcal D_{y}^{\Sigma_{yy}} \! = \mathcal I^N_x \otimes \mathcal D^N_y,
\end{array}
\end{equation}
where $\mathcal I$ symbolizes a 1D identity matrix.
In \eqref{Tensor_product_norm_matrices} and \eqref{Tensor_product_difference_matrices}, superscript $^N$ indicates the `{\it reference}' grid that aligns with the boundaries while $^M$ indicates the grid that is staggered with respect to the `{\it reference}' grid. 

In the $x$ direction, the 1D norm matrices, i.e., $\mathcal A^M_x$ and $\mathcal A^N_x$, are simply identity matrices scaled by the grid spacing $\Delta x$, while the 1D finite difference matrices, i.e., $\mathcal D^M_x$ and $\mathcal D^N_x$, are built from the standard stencil $\nicefrac{ \left[ \nicefrac{1}{24}, \ -\nicefrac{9}{8}, \ \nicefrac{9}{8}, \ -\nicefrac{1}{24} \right] }{ \Delta x }$, which is `{\it wrapped around}' when approaching boundaries to account for periodic boundary condition.
Constructed as such, these operators satisfy the following relation:
\begin{equation}
\label{Operator_relation_1D_x}
\small
\mathcal A^N_x \mathcal D^M_x + \left( \mathcal A^M_x \mathcal D^N_x\right)^T = \bf 0.
\end{equation}

In the $y$ direction, we use the 1D operators provided in \cite[p.6]{gao2019sbp}\footnotemark[2] for $\mathcal A^M_y$, $\mathcal A^N_y$, $\mathcal D^M_y$, and $\mathcal D^N_y$, which satisfy the following relation:
\begin{equation}
\small
\label{Operator_relation_1D_y}
\mathcal A^N_y \mathcal D^M_y + \left( \mathcal A^M_y \mathcal D^N_y\right)^T 
=
\mathcal E_y^R \left(\mathcal P_y^R\right)^T - \ \mathcal E_y^L \left(\mathcal P_y^L\right)^T,
\end{equation}
where $\mathcal E_y^R$ and $\mathcal E_y^L$ are canonical basis vectors that select values of solution variables defined on the $N$ grid at the top and bottom boundaries, respectively; $\mathcal P_y^R$ and $\mathcal P_y^L$ are projection vectors that extrapolate values of solution variables defined on the $M$ grid to the top and bottom boundaries, respectively.
\footnotetext[2]{We use this particular set of operators to demonstrate the methodology of interface treatment. Other alternative choices for these operators exist.}

The discrete energy associated with semi-discretized system \eqref{Elastic_wave_equation_Discretization_2D_analysis} is defined as:
\begin{equation}
\label{Discrete_energy}
\small
E = \tfrac{1}{2} V_i^T \left( \mathcal A^{V_i}  \boldsymbol{\rho}^{V_i} \right) V_i^{\phantom{T}} 
+
\tfrac{1}{2} \Sigma_{ij}^T \left( \mathcal A^{\Sigma_{ij}}  S^{\Sigma_{kl}}_{ijkl} \right) \Sigma_{kl}^{\phantom{T}}\ ,
\end{equation}
where the index summation convention applies only to $i$, $j$, $k$, and $l$ in the subscripts.
Differentiating $E$ with respect to time $t$ and substituting the equations in \eqref{Elastic_wave_equation_Discretization_2D_analysis}, it can be shown that 
\begin{equation}
\label{Discrete_energy_analysis}
\small
\begin{array}{rclll}
\displaystyle \frac{d E}{d t} \!
&=& \!
\displaystyle V_x^T \left[ \mathcal I^M_x \! \otimes \mathcal E^R_y \right] \mathcal A^M_x \left[ \mathcal I^M_x \! \otimes ( \mathcal P^R_y )^T \right] \Sigma_{xy}^{\phantom{T}} \!
&+& 
\displaystyle \Sigma_{yy}^T \left[ \mathcal I^N_x \! \otimes \mathcal E^R_y \right] \mathcal A^N_x \left[ \mathcal I^N_x \! \otimes ( \mathcal P^R_y )^T \right] V_y^{\phantom{T}} \\[2.ex]
&-& \!
\displaystyle V_x^T \left[ \mathcal I^M_x \! \otimes \mathcal E^L_y \right] \mathcal A^M_x \left[ \mathcal I^M_x \! \otimes ( \mathcal P^L_y )^T \right] \Sigma_{xy}^{\phantom{T}} \!
&-&
\displaystyle \Sigma_{yy}^T \left[ \mathcal I^N_x \! \otimes \mathcal E^L_y \right] \mathcal A^N_x \left[ \mathcal I^N_x \! \otimes ( \mathcal P^L_y )^T \right] V_y^{\phantom{T}}. 
\end{array}
\end{equation} 
With the above discrete energy analysis result, we can now modify system \eqref{Elastic_wave_equation_Discretization_2D_analysis} accordingly to account for boundary and interface conditions.

In the following, we use superscripts $^+$ and $^-$ to distinguish systems or terms from the upper and lower regions of Figure \ref{Grid_Configuration_Interface}, respectively.
To account for the free-surface boundary condition on the top boundary, i.e., $\sigma_{xy} = \sigma_{yy} = 0$, the first two equations in the upper region system are appended with penalty terms, i.e., SATs, as follows:
\begin{equation}
\label{SATs_free_surface_upper}
\small
\text{(+)}
\left\{
\begin{array}{rcl}
\displaystyle \mathcal A^{V_x} \boldsymbol{\rho}^{V_x} \frac{d V_x}{d t} &=& \displaystyle \mathcal A^{V_x} \mathcal D^{\Sigma_{xx}}_x \Sigma_{xx} + \mathcal A^{V_x} \mathcal D^{\Sigma_{xy}}_y \Sigma_{xy} \\
&+& \displaystyle \eta_T^{V_x} \left[ \mathcal I^M_x \otimes \mathcal E^R_y \right] \mathcal A^M_x \left\{ \left[ \mathcal I^M_x \otimes ( \mathcal P^R_y )^T \right] \Sigma_{xy} - \boldsymbol 0 \right\}; \\[1.ex]
\displaystyle \mathcal A^{V_y} \boldsymbol{\rho}^{V_y} \frac{d V_y}{d t} &=& \displaystyle \mathcal A^{V_y} \mathcal D^{\Sigma_{xy}}_x \Sigma_{xy} + \mathcal A^{V_y} \mathcal D^{\Sigma_{yy}}_y \Sigma_{yy} \\
&+& \eta_T^{V_y} \left[ \mathcal I^N_x \otimes \mathcal P^R_y \right] \mathcal A^N_x \displaystyle \left\{ \left[ \mathcal I^N_x \otimes ( \mathcal E^R_y )^T \right] \Sigma_{yy} - \boldsymbol 0 \right\},
\end{array}
\right.
\end{equation}
where $\eta_T^{V_x} = \eta_T^{V_y} = -1$ are penalty parameters. Similarly, the first two equations of the lower region system are modified as follows to account for the free-surface boundary condition on the bottom boundary:
\begin{equation}
\label{SATs_free_surface_lower}
\small
\text{(-)}
\left\{
\begin{array}{rcl}
\displaystyle \mathcal A^{V_x} \boldsymbol{\rho}^{V_x} \frac{d V_x}{d t} &=& \displaystyle \mathcal A^{V_x} \mathcal D^{\Sigma_{xx}}_x \Sigma_{xx} + \mathcal A^{V_x} \mathcal D^{\Sigma_{xy}}_y \Sigma_{xy} \\
&+& \displaystyle \eta_B^{V_x} \left[ \mathcal I^M_x \otimes \mathcal E^L_y \right] \mathcal A^M_x \left\{ \left[ \mathcal I^M_x \otimes ( \mathcal P^L_y )^T \right] \Sigma_{xy} - \boldsymbol 0 \right\}; \\[1.ex]
\displaystyle \mathcal A^{V_y} \boldsymbol{\rho}^{V_y} \frac{d V_y}{d t} &=& \displaystyle \mathcal A^{V_y} \mathcal D^{\Sigma_{xy}}_x \Sigma_{xy} + \mathcal A^{V_y} \mathcal D^{\Sigma_{yy}}_y \Sigma_{yy} \\
&+& \displaystyle \eta_B^{V_y} \left[ \mathcal I^N_x \otimes \mathcal P^L_y \right] \mathcal A^N_x \left\{ \left[ \mathcal I^N_x \otimes ( \mathcal E^L_y )^T \right] \Sigma_{yy} - \boldsymbol 0 \right\},
\end{array}
\right.
\end{equation}
where $\eta_B^{V_x} = \eta_B^{V_y} = 1$ are penalty parameters. 

\begin{figure}[h]
\vspace{-.75em}
\captionsetup{width=1\textwidth, font=small,labelfont=small}
\centering
\begin{subfigure}[b]{0.495\textwidth}
\centering\includegraphics[scale=0.125]{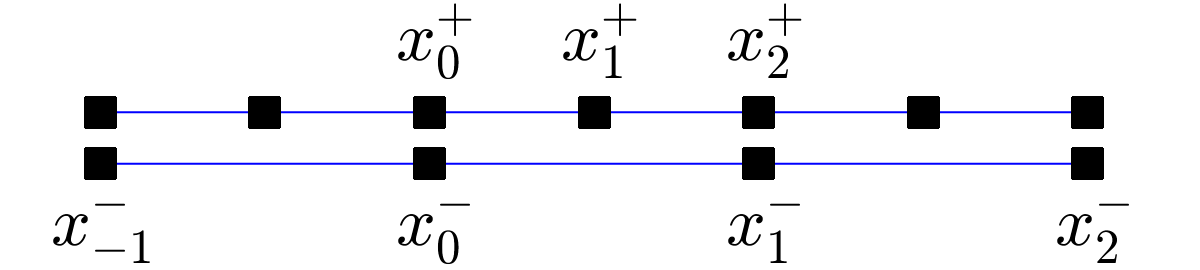}
\vspace{-0.25em}
\end{subfigure}
\begin{subfigure}[b]{0.495\textwidth}
\centering\includegraphics[scale=0.125]{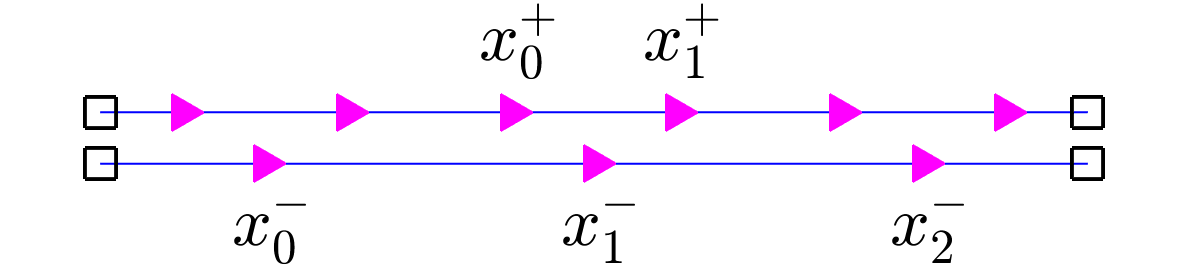}
\vspace{-0.25em}
\end{subfigure} 
\caption{Illustration of grid points involved for transfer operators $\mathcal{T}^{N_{-}^{+}}$ (left) and $\mathcal{T}^{M_{-}^{+}}$ (right).} 
\label{Grid_configuration_transfer_operator}
\vspace{-1.5em}
\end{figure}

To account for the interface conditions, i.e., $\sigma^{+}_{xy} \!=\! \sigma^{-}_{xy}$, $\sigma^{+}_{yy} \!=\! \sigma^{-}_{yy}$, $v^{+}_{x} \!=\! v^{-}_{x}$, and $v^{+}_{y} \!=\! v^{-}_{y}$ (see \cite[p.52]{stein2009introduction}) the upper region system is further modified by appending additional SATs as follows:
\begin{equation}
\label{SATs_interface_upper}
\small
\text{(+)}
\left\{
\begin{array}{l}
\displaystyle \mathcal A^{V^+_x} \boldsymbol{\rho}^{V^+_x} \frac{d V^+_x}{d t} = \displaystyle \mathcal A^{V^+_x} \mathcal D^{\Sigma^+_{xx}}_x \Sigma^+_{xx} + \mathcal A^{V^+_x} \mathcal D^{\Sigma^+_{xy}}_y \Sigma^+_{xy} \\[1.5ex]
\enskip + \ \displaystyle \eta_I^{V^+_x} \left[ \mathcal I^{M^+}_x \!\! \otimes \mathcal E^{L^+}_y \right] \mathcal A^{M^+}_x \!\! \left\{ \left[ \mathcal I^{M^+}_x \!\! \otimes ( \mathcal P^{L^+}_y )^T \right] \Sigma^+_{xy} - \mathcal{T}^{M_{-}^{+}} \left( \left[ \mathcal I^{M^-}_x \!\!\otimes ( \mathcal P^{R^-}_y )^T \right] \Sigma^-_{xy} \right) \right\}; \\[2ex]
\displaystyle \mathcal A^{V^+_y} \boldsymbol{\rho}^{V^+_y} \frac{d V^+_y}{d t} = \displaystyle \mathcal A^{V^+_y} \mathcal D^{\Sigma^+_{xy}}_x \Sigma^+_{xy} + \mathcal A^{V^+_y} \mathcal D^{\Sigma^+_{yy}}_y \Sigma^+_{yy} \\[1.5ex]
\enskip + \ \displaystyle \eta_I^{V^+_y} \left[ \mathcal I^{N^+}_x \!\! \otimes \mathcal P^{L^+}_y \right]  \mathcal A^{N^+}_x \!\! \left\{ \left[ \mathcal I^{N^+}_x \!\! \otimes ( \mathcal E^{L^+}_y )^T \right] \Sigma^+_{yy} - \mathcal{T}^{N_{-}^{+}} \left( \left[ \mathcal I^{N^-}_x \!\! \otimes ( \mathcal E^{R^-}_y )^T \right] \Sigma^-_{yy} \right) \right\}; \\[2ex]
\displaystyle \mathcal A^{\Sigma^+_{xx}} S^{\Sigma^+_{kl}}_{xxkl} \frac{d \Sigma^+_{kl}}{d t} = \displaystyle \mathcal A^{\Sigma^+_{xx}} \mathcal D^{V^+_x}_x V^+_x; \\[2ex]
\displaystyle \mathcal A^{\Sigma^+_{xy}} S^{\Sigma^+_{kl}}_{xykl} \frac{d \Sigma^+_{kl}}{d t} = \displaystyle \tfrac{1}{2} \mathcal A^{\Sigma^+_{xy}} \left( \mathcal D^{V^+_y}_x V^+_y + \mathcal D^{V^+_x}_y V^+_x \right) \\[1.5ex]
\enskip + \ \displaystyle \tfrac{1}{2} \eta_I^{\Sigma^+_{xy}} \left[ \mathcal I^{M^+}_x \!\! \otimes \mathcal P^{L^+}_y \right] \mathcal A^{M^+}_x \!\! \left\{ \left[ \mathcal I^{M^+}_x \!\! \otimes ( \mathcal E^{L^+}_y )^T \right] V^+_x - \mathcal{T}^{M_{-}^{+}} \left( \left[ \mathcal I^{M^-}_x \!\! \otimes ( \mathcal E^{R^-}_y )^T \right] V^-_x \right) \right\}; \\[2ex]
\displaystyle \mathcal A^{\Sigma^+_{yy}} S^{\Sigma^+_{kl}}_{yykl} \frac{d \Sigma^+_{kl}}{d t} = \displaystyle \mathcal A^{\Sigma^+_{yy}} \mathcal D^{V^+_y}_y V^+_y \\[1.5ex]
\enskip + \ \displaystyle \eta_I^{\Sigma^+_{yy}} \left[ \mathcal I^{N^+}_x \!\! \otimes \mathcal E^{L^+}_y \right] \mathcal A^{N^+}_x \!\! \left\{ \left[ \mathcal I^{N^+}_x \!\! \otimes ( \mathcal P^{L^+}_y )^T \right] V^+_y - \mathcal{T}^{N_{-}^{+}} \left( \left[ \mathcal I^{N^-}_x \!\! \otimes ( \mathcal P^{R^-}_y )^T \right] V^-_y \right) \right\},
\end{array}
\right.
\end{equation}
while the lower region system is further modified by appending additional SATs as follows:
\begin{equation}
\label{SATs_interface_lower}
\small
\text{(-)}
\left\{
\begin{array}{l}
\displaystyle \mathcal A^{V^-_x} \boldsymbol{\rho}^{V^-_x} \frac{d V^-_x}{d t} = \displaystyle \mathcal A^{V^-_x} \mathcal D^{\Sigma^-_{xx}}_x \Sigma^-_{xx} + \mathcal A^{V^-_x} \mathcal D^{\Sigma^-_{xy}}_y \Sigma^-_{xy} \\[1.5ex]
\enskip + \ \displaystyle \eta_I^{V^-_x} \left[ \mathcal I^{M^-}_x \!\! \otimes \mathcal E^{R^-}_y \right] \mathcal A^{M^-}_x \!\! \left\{ \left[ \mathcal I^{M^-}_x \!\! \otimes ( \mathcal P^{R^-}_y )^T \right] \Sigma^-_{xy} - \mathcal{T}^{M_{+}^{-}} \left( \left[ \mathcal I^{M^+}_x \!\!\otimes ( \mathcal P^{L^+}_y )^T \right] \Sigma^+_{xy} \right) \right\}; \\[2ex]
\displaystyle \mathcal A^{V^-_y} \boldsymbol{\rho}^{V^-_y} \frac{d V^-_y}{d t} = \displaystyle \mathcal A^{V^-_y} \mathcal D^{\Sigma^-_{xy}}_x \Sigma^-_{xy} + \mathcal A^{V^-_y} \mathcal D^{\Sigma^-_{yy}}_y \Sigma^-_{yy} \\[1.5ex]
\enskip + \ \displaystyle \eta_I^{V^-_y} \left[ \mathcal I^{N^-}_x \!\! \otimes \mathcal P^{R^-}_y \right]  \mathcal A^{N^-}_x \!\! \left\{ \left[ \mathcal I^{N^-}_x \!\! \otimes ( \mathcal E^{R^-}_y )^T \right] \Sigma^-_{yy} - \mathcal{T}^{N_{+}^{-}} \left( \left[ \mathcal I^{N^+}_x \!\! \otimes ( \mathcal E^{L^+}_y )^T \right] \Sigma^+_{yy} \right) \right\}; \\[2ex]
\displaystyle \mathcal A^{\Sigma^-_{xx}} S^{\Sigma^-_{kl}}_{xxkl} \frac{d \Sigma^-_{kl}}{d t} = \displaystyle \mathcal A^{\Sigma^-_{xx}} \mathcal D^{V^-_x}_x V^-_x; \\[2ex]
\displaystyle \mathcal A^{\Sigma^-_{xy}} S^{\Sigma^-_{kl}}_{xykl} \frac{d \Sigma^-_{kl}}{d t} = \displaystyle \tfrac{1}{2} \mathcal A^{\Sigma^-_{xy}} \left( \mathcal D^{V^-_y}_x V^-_y + \mathcal D^{V^-_x}_y V^-_x \right) \\[1.5ex]
\enskip + \ \displaystyle \tfrac{1}{2} \eta_I^{\Sigma^-_{xy}} \left[ \mathcal I^{M^-}_x \!\! \otimes \mathcal P^{R^-}_y \right] \mathcal A^{M^-}_x \!\! \left\{ \left[ \mathcal I^{M^-}_x \!\! \otimes ( \mathcal E^{R^-}_y )^T \right] V^-_x - \mathcal{T}^{M_{+}^{-}} \left( \left[ \mathcal I^{M^+}_x \!\! \otimes ( \mathcal E^{L^+}_y )^T \right] V^+_x \right) \right\}; \\[2ex]
\displaystyle \mathcal A^{\Sigma^-_{yy}} S^{\Sigma^-_{kl}}_{yykl} \frac{d \Sigma^-_{kl}}{d t} = \displaystyle \mathcal A^{\Sigma^-_{yy}} \mathcal D^{V^-_y}_y V^-_y \\[1.5ex]
\enskip + \ \displaystyle \eta_I^{\Sigma^-_{yy}} \left[ \mathcal I^{N^-}_x \!\! \otimes \mathcal E^{R^-}_y \right] \mathcal A^{N^-}_x \!\! \left\{ \left[ \mathcal I^{N^-}_x \!\! \otimes ( \mathcal P^{R^-}_y )^T \right] V^-_y - \mathcal{T}^{N_{+}^{-}} \left( \left[ \mathcal I^{N^+}_x \!\! \otimes ( \mathcal P^{L^+}_y )^T \right] V^+_y \right) \right\},
\end{array}
\right.
\end{equation}
where {\small $\eta_L^{V^+_x} \! = \eta_L^{\Sigma^+_{xy}} \! = \eta_L^{V^+_y} \! = \eta_L^{\Sigma^+_{yy}} \! = \tfrac{1}{2}$} and {\small $\eta_R^{V^-_x} \! = \eta_R^{\Sigma^-_{xy}} \! = \eta_R^{V^-_y} \! = \eta_R^{\Sigma^-_{yy}} \! = - \tfrac{1}{2}$} are proper choices for these penalty parameters. 
Moreover, {\small $\mathcal{T}^{M_{-}^{+}}$}, {\small $\mathcal{T}^{N_{-}^{+}}$} in \eqref{SATs_interface_upper} and {\small $\mathcal{T}^{M_{+}^{-}}$}, {\small $\mathcal{T}^{N_{+}^{-}}$} in \eqref{SATs_interface_lower} are transfer operators that satisfy the following relations:
\begin{equation}
\small
\label{Relations_interpolation_operator}
\mathcal A^{N^+}_x \mathcal{T}^{N_{-}^{+}} = \left( \mathcal A^{N^-}_x \mathcal{T}^{N_{+}^{-}} \right)^T \text{ and } \,\, 
\mathcal A^{M^+}_x \mathcal{T}^{M_{-}^{+}} = \left( \mathcal A^{M^-}_x \mathcal{T}^{M_{+}^{-}} \right)^T\!. 
\end{equation}
These transfer operators operate on the interface only, e.g., $\mathcal{T}^{N_{-}^{+}}$ interpolates from lower region $N$ grid points on the interface to upper region $N$ grid points on the interface.  
They are usually designed with the assistance of symbolic computing softwares. 
For the interface illustrated in Figure \ref{Grid_Configuration_Interface}, which has a contrast ratio $1 \!\! : \!\! 2$ in grid spacing, the operators $\mathcal{T}^{N_{-}^{+}}$ and $\mathcal{T}^{M_{-}^{+}}$ that we use can be characterized by the formulas in \eqref{Transfer_operator_N} and \eqref{Transfer_operator_M}, respectively, for the small collections of grid points illustrated in Figure \ref{Grid_configuration_transfer_operator}, thanks to the repeated grid patterns. Moreover, $\mathcal{T}^{N_{+}^{-}}$ and $\mathcal{T}^{M_{+}^{-}}$ can be derived from $\mathcal{T}^{N_{-}^{+}}$ and $\mathcal{T}^{M_{-}^{+}}$, respectively, via the relations in \eqref{Relations_interpolation_operator}.
As in the case of SBP operators, these transfer operators are not unique, either. 
\begin{equation}
\label{Transfer_operator_N}
\small
f(x^+_0) \leftarrow f(x^-_0); 
\quad 
f(x^+_1) \leftarrow -\tfrac{1}{16} f(x^-_{-1}) 
			     +\tfrac{9}{16} f(x^-_0) 
			     +\tfrac{9}{16} f(x^-_1) 
			     -\tfrac{1}{16} f(x^-_2);
\quad 
f(x^+_2) \leftarrow f(x^-_1).
\end{equation} 
\\[-5ex]
\begin{equation}
\label{Transfer_operator_M}
\small
f(x^+_0) \leftarrow  \tfrac{5}{32} f(x^-_0) 
			    +\tfrac{15}{16} f(x^-_1) 
			     -\tfrac{3}{32} f(x^-_2); 
\qquad \qquad
f(x^+_1) \leftarrow  -\tfrac{3}{32} f(x^-_0) 
			     +\tfrac{15}{16} f(x^-_1) 
			     +\tfrac{5}{32} f(x^-_2).
\end{equation}
\vspace{0.25em}

With the above choices on the SATs, it can be verified that the overall semi-discretization conserves the discrete energy $E$ defined in \eqref{Discrete_energy}.
Now that the proper SATs have been derived, we can remove the norm matrices by dividing them from both sides of the equations; see (\ref{SATs_free_surface_upper}-\ref{SATs_interface_lower}).
From the implementation perspective, the appended SATs amount to modifying the corresponding derivative approximations, e.g., the SAT in the first equation of \eqref{SATs_interface_lower} modifies {\small $\mathcal D^{\Sigma^-_{xy}}_y \Sigma^-_{xy}$}.
With this understanding, the above semi-discretizations for system \eqref{Elastic_wave_equation_PDE_2D_analysis} can easily be reverted to forms that conform to system \eqref{Elastic_wave_equation_PDE_2D_implementation} by first absorbing the appended SATs with modified derivative approximations.

\section{Numerical examples} \label{section_examples}
The first numerical example concerns a homogeneous medium characterized by parameters $\rho = 1\, \text{kg}/\text{m}^3$, $c_p = 2\, \text{m}/\text{s}$ and $c_s = 1\, \text{m}/\text{s}$.
The grid spacings of the upper and lower regions are chosen as $0.004\, \text{m}$ and $0.008\, \text{m}$, respectively, while the time step length is chosen as $0.001\, \text{s}$.
The rest of the numerical setup is the same as for the first example of \cite{gao2019sbp}, including sizes of the grids, source and receiver locations, as well as the source profile.
\begin{figure}[h]
\vspace{-0.75em}
\captionsetup{width=1\textwidth, font=small,labelfont=small}
\centering
\begin{subfigure}[b]{0.495\textwidth}
\centering\includegraphics[scale=0.1]{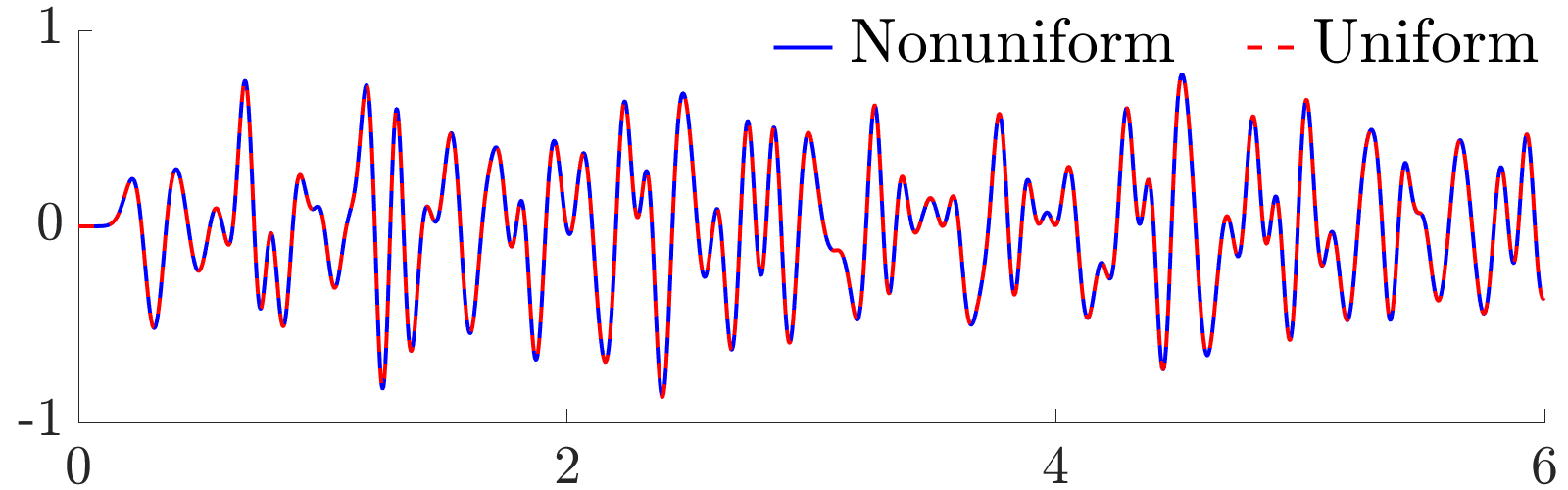}
\end{subfigure}
\begin{subfigure}[b]{0.495\textwidth}
\centering\includegraphics[scale=0.1]{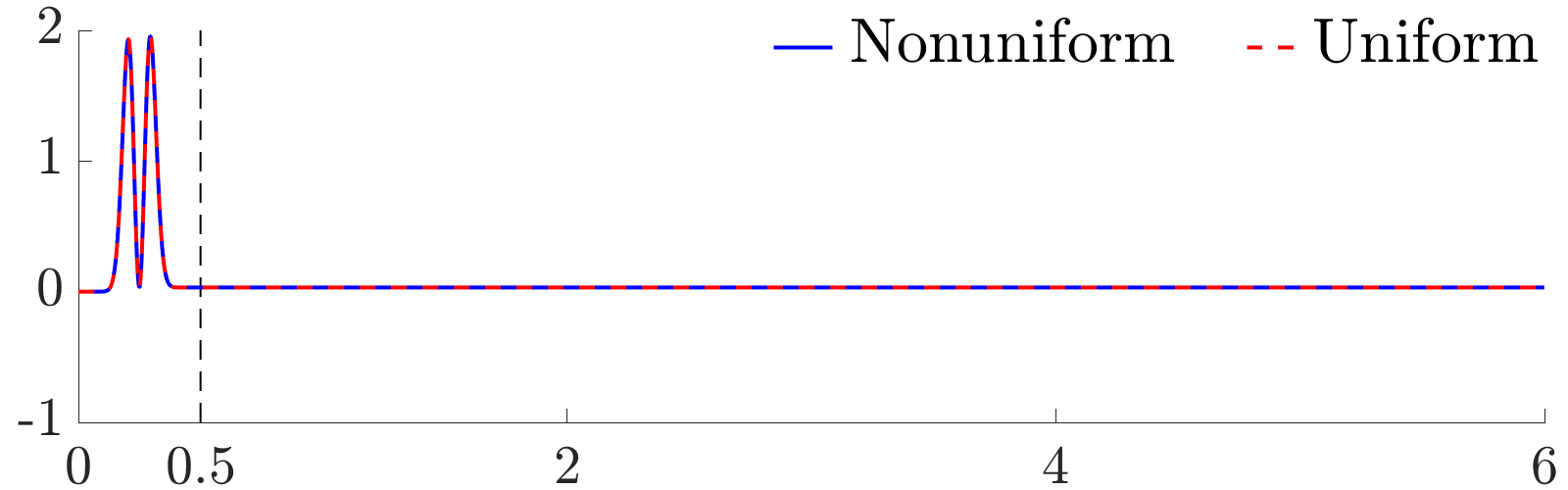}
\end{subfigure} 
\caption{Seismogram (left) and evolvement of discrete energy (right); Homogeneous media.} 
\label{Homogeneous_seismogram_and_energy}
\vspace{-1.em}
\end{figure}

The recorded seismogram and evolvement of discrete energy for the first $6\, \text{s}$ are displayed in Figure \ref{Homogeneous_seismogram_and_energy}, where we observe good agreement between the uniform grid simulation result and nonuniform grid simulation result using the presented SBP-SAT approach.
The source term $\mathcal S$ in \eqref{Elastic_wave_equation_PDE_2D_implementation}, which is omitted from the analysis, is responsible for the initial `{\it bumps}' in the evolvement of discrete energy.
After the source effect tapers off (at around $0.5\, \text{s}$), the discrete energy remains constant as expected (with a value at approximately $0.0318$).

\begin{figure}[h]
\captionsetup{width=1\textwidth, font=small,labelfont=small}
\centering
\begin{subfigure}[b]{0.495\textwidth}
\centering\includegraphics[scale=0.0475]{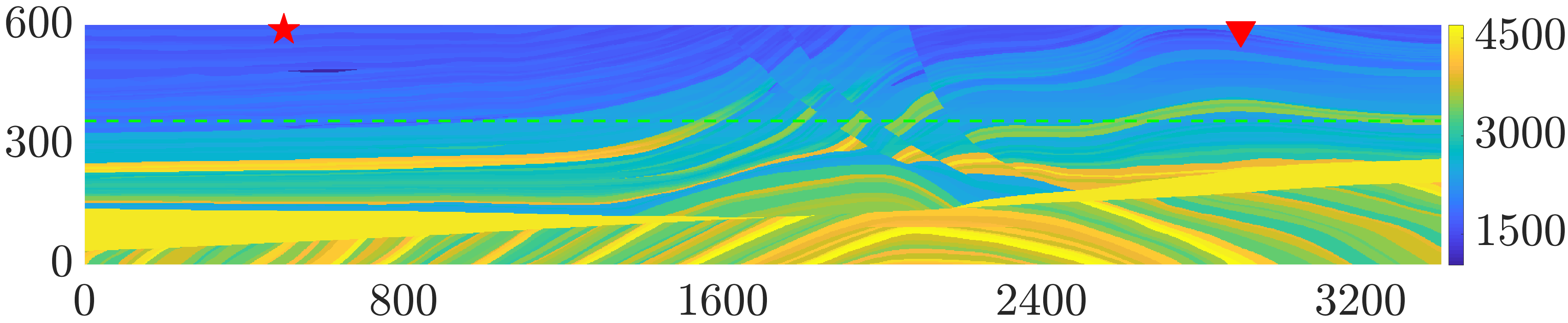}
\end{subfigure} 
\begin{subfigure}[b]{0.495\textwidth}
\centering\includegraphics[scale=0.0475]{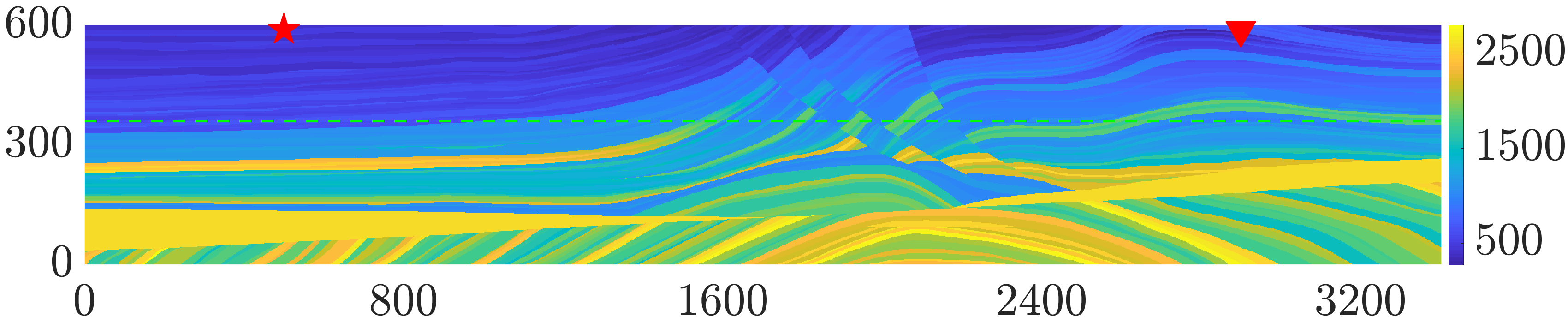}
\end{subfigure} 
\caption{Media parameters $c_p$ (left) and $c_s$ (right); Colorbar reflects wave-speed with unit m/s.}
\label{Marmousi2_parameter}
\vspace{1em}
\end{figure}

\begin{figure}[h]
\captionsetup{width=1\textwidth, font=small,labelfont=small}
\centering
\begin{subfigure}[b]{0.495\textwidth}
\centering\includegraphics[scale=0.1]{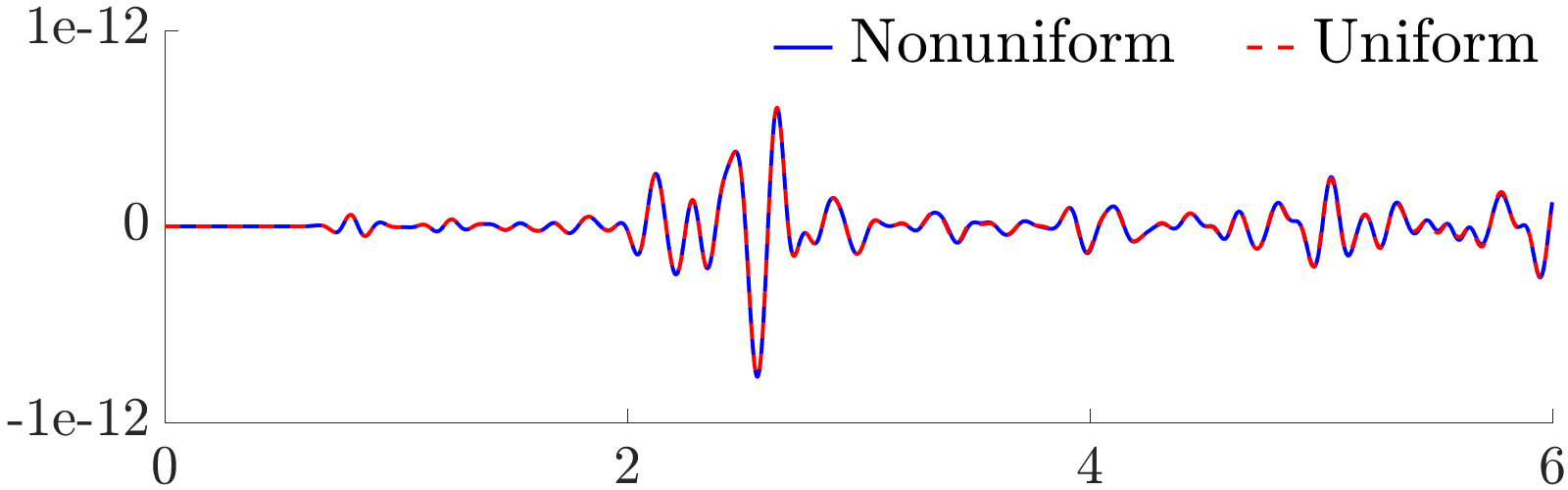}
\end{subfigure}
\begin{subfigure}[b]{0.495\textwidth}
\centering\includegraphics[scale=0.1]{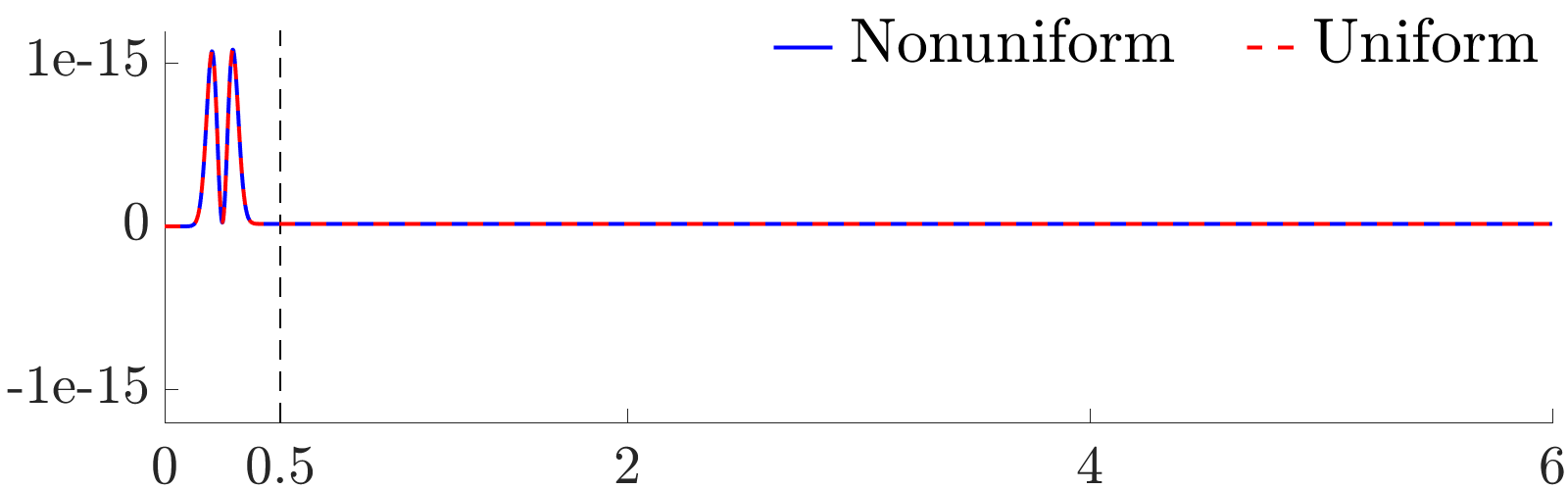}
\end{subfigure} 
\caption{Seismogram (left) and evolvement of discrete energy (right); Heterogeneous media.}
\label{Heterogeneous_seismogram_and_energy}
\end{figure}

The second numerical example concerns a heterogeneous medium downsampled from the Marmousi2 model; see \cite{martin2006marmousi2}.
Media parameters $c_p$ and $c_s$ used here are illustrated in Figure \ref{Marmousi2_parameter}.\footnotemark[3]
Grid spacing is chosen as $2\, \text{m}$ and $4\, \text{m}$ for upper and lower regions, respectively. 
Time step length is chosen as $2$e-$4\, \text{s}$ and $3$e-$4\, \text{s}$ for uniform and nonuniform grid simulations, respectively.
Same plots as in the previous example are displayed in Figure \ref{Heterogeneous_seismogram_and_energy}, from where similar observations can be made.
\footnotetext[3]{For the purpose of illustration, the distance between neighboring parameter grid points is assigned a value $2\, \text{m}$, which is the same as the grid spacing used in uniform grid simulation. 
Bilinear interpolation is used when discretization grid points do not match parameter grid points due to grid staggering.
We note here that media parameters for uniform grid and nonuniform grid simulations are sampled differently; thus, small discrepancies in simulation results should be allowed.}

\section{Summary}\label{section_summary}
Finite difference discretization of isotropic elastic wave equations is considered. 
An interface treatment procedure is presented to connect two uniformly discretized regions of different grid spacings, where interface conditions are imposed through carefully designed penalty terms, which are also referred to as the simultaneous approximation terms.
The overall semi-discretization conserves a discrete energy resembling the continuous energy possessed by the elastic wave system, which is demonstrated on both homogeneous and heterogeneous media. 


\bibliographystyle{unsrtnat}
\bibliography{References}

\end{document}